\documentclass{svmult}
\usepackage{newtxtext}
\usepackage{newtxmath}
\usepackage{makeidx}         
\usepackage{graphicx}        

\usepackage{multicol}        
\usepackage[bottom]{footmisc}
\usepackage{csml}

\usepackage[pdftex,colorlinks=true,allcolors=blue]{hyperref}

\makeindex
\begin{document}
%
%
%
%
\title*{Manifold-based B-splines on unstructured meshes }
\titlerunning{Manifold-based B-splines}  
%
\author{Qiaoling Zhang, Thomas Takacs and Fehmi Cirak}
\authorrunning{Zhang, Q., Takacs, T. and Cirak, F} 
%
\tocauthor{Qiaoling Zhang, Thomas Takacs and Fehmi Cirak}
\institute{Qiaoling Zhang, Fehmi Cirak \at Department of Engineering, University of Cambridge, Cambridge, CB2 1PZ, UK,  \\ \mbox{\email{zq217@cam.ac.uk}}, \mbox{\email{fc286@cam.ac.uk}} \and
Thomas Takacs \at Institute of Applied Geometry, Johannes Kepler University, Linz, A-4040, Austria, \\ \mbox{\email{thomas.takacs@jku.at}} }

\maketitle              

\abstract{
We introduce new manifold-based splines that are able to exactly reproduce B-splines on unstructured surface meshes. Such splines can be used in isogeometric analysis (IGA) to represent smooth surfaces of arbitrary topology. Since prevalent computer-aided design (CAD) models are composed of tensor-product B-spline patches, any IGA suitable construction should be able to reproduce B-splines. 
To achieve this goal, we focus on univariate manifold-based constructions that can reproduce B-splines. 
The manifold-based splines are constructed by smoothly blending together polynomial interpolants defined on overlapping charts. The proposed constructions automatically reproduce B-splines in regular parts of the mesh, with no extraordinary vertices, and polynomial basis functions in the remaining parts of the mesh. We study and compare analytically and numerically the finite element convergence of several univariate constructions. The obtained results directly carry over to the tensor-product case.}

\newpage
%
%
\section{Introduction}
%
Manifold-based surface construction techniques from geometric modelling provide an elegant and flexible framework for generating basis functions on surfaces with arbitrary topology~\cite{grimm1995modeling,navau2000modeling,ying2004simple,della2008construction}. They combine manifold descriptions from differential geometry, see e.g.~\cite{schutz1980geometrical,doCarmo1976differential},  with the flexibility of the partition of unity framework from numerical analysis~\cite{Melenk1996a,babuvska1997partition}. If a manifold surface in~$\mathbb R^3$ can be mapped onto a single planar parametric domain in~$\mathbb R^2$, it is straightforward to obtain partition of unity basis functions of any desired regularity on the parametric domain. Although it is impossible to map a surface with arbitrary topology onto a single parametric domain, it can always be represented as an atlas composed of a number of charts. The charts consist of a planar domains in~$\mathbb R^2$ that map onto the manifold surface in~$\mathbb R^3$. The planar chart domains in~$\mathbb R^2$ are not connected and transition functions are used to navigate between the different domains. The manifold-based basis functions are obtained by simply applying the partition of unity method on the collection of chart domains~\cite{majeedCirak:2016, zhang:2019}. The flexibility of the original partition of unity method carries over to the manifold case.  The partition of unity functions, referred to as blending functions, in this paper, the local approximants on each chart domain and the transition functions can all be chosen to fit the requirements of the application at hand. 

\begin{figure}[]
    \centering
	 \includegraphics[width=0.55\textwidth]{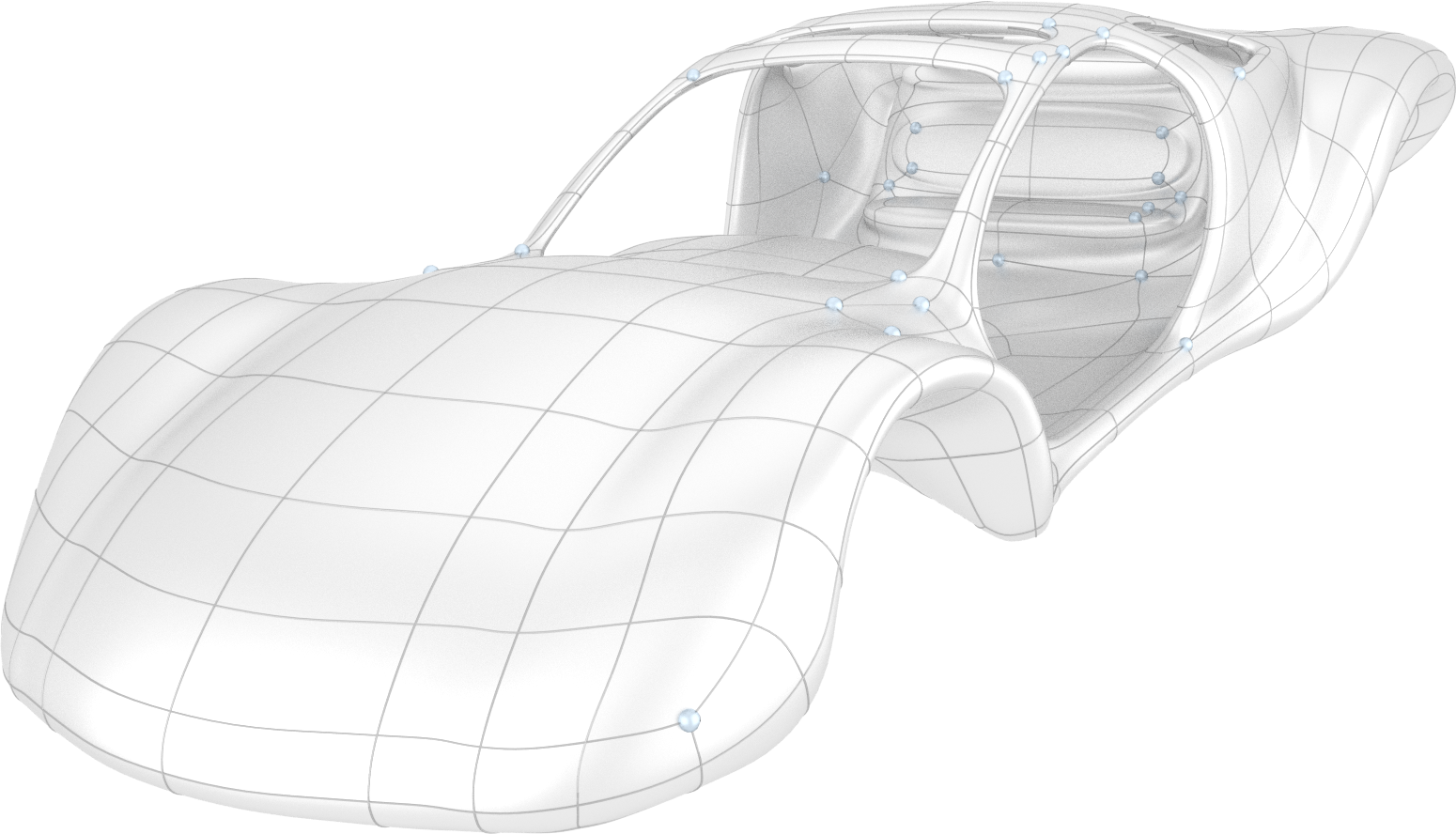} \\ [1em]
	  \includegraphics[width=0.45\textwidth]{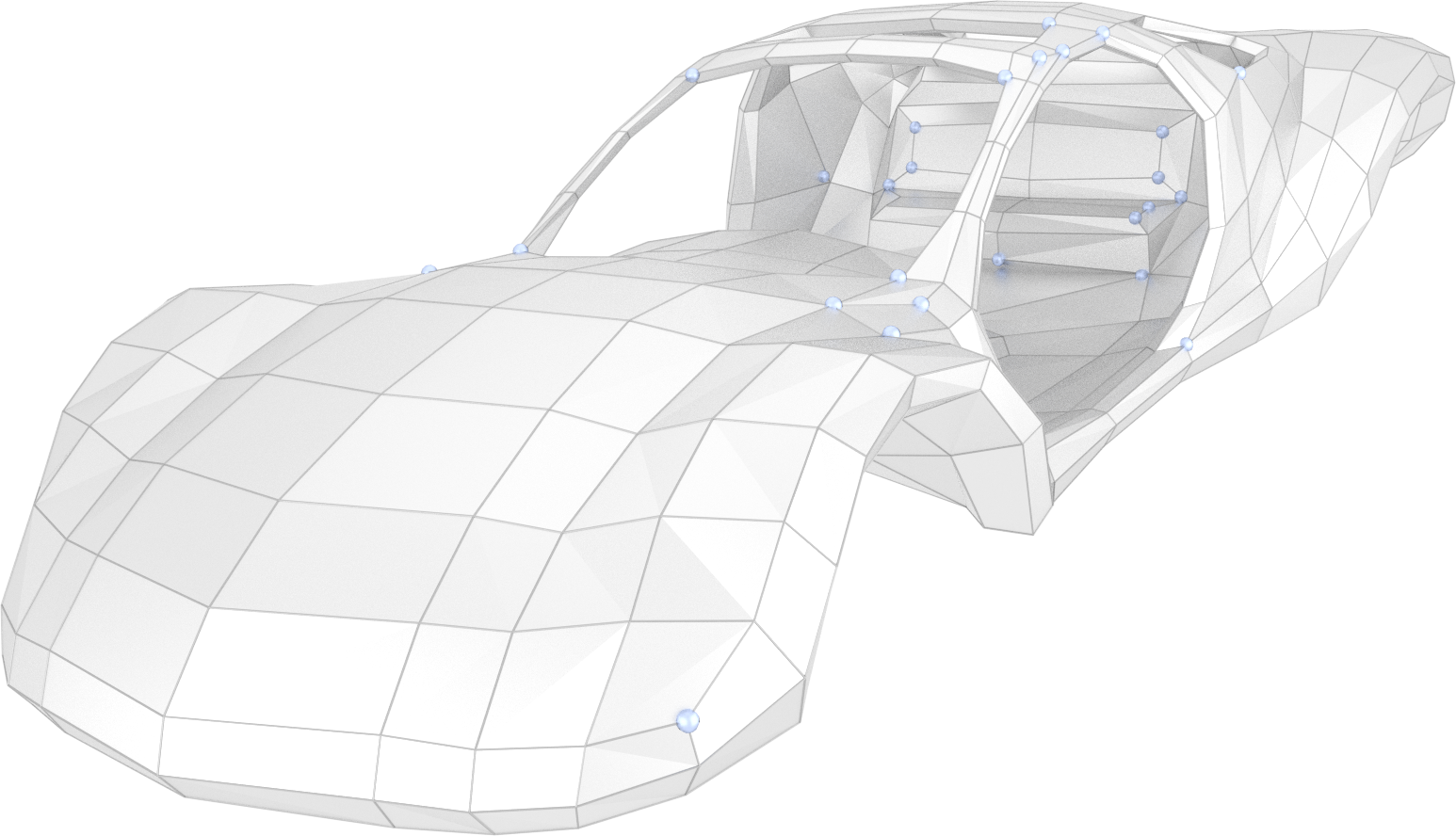} \hfill \includegraphics[width=0.45\textwidth]{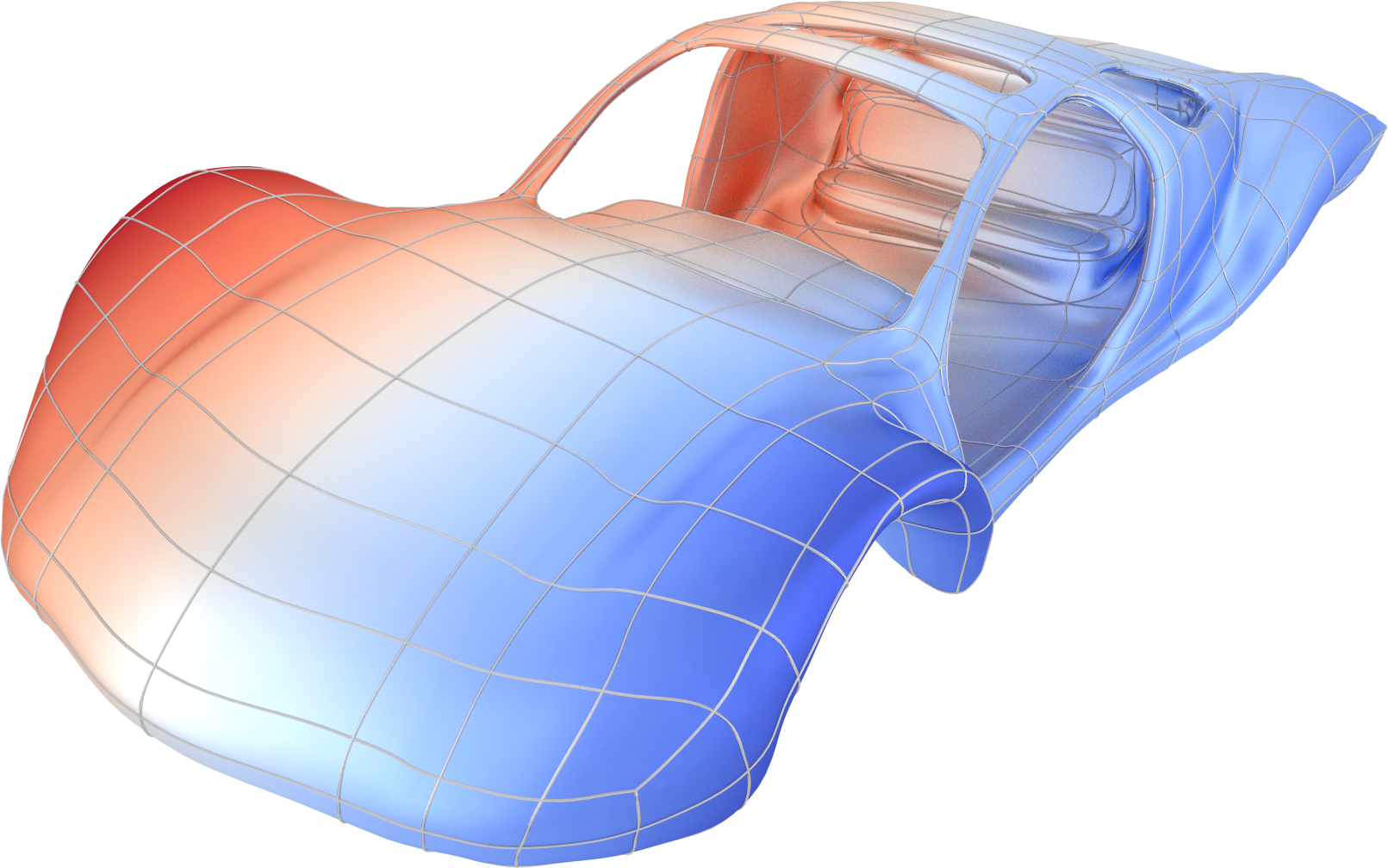} 
 	\caption{Isogeometric analysis of a car body as used in computer animation. The manifold-based representation (top) is obtained from an unstructured coarse quadrilateral control mesh (bottom, left), with extraordinary vertices indicated by blue spheres. The deflected shape of the car body subjected to an axial torsion is computed with Kirchhoff-Love shell finite elements (bottom, right). \label{fig:car}}
\end{figure}

The definition of smooth functions over unstructured meshes, such as shown in Figure~\ref{fig:car}, or multi-patch geometries has always been of vital interest in the context of isogeometric analysis, cf.~\cite{Hughes:2005aa,cottrell2009isogeometric}. Approaches for the definition of such functions include the subdivision surfaces~\cite{Cirak:2000aa, Cirak:2011aa,Peters:2008aa,zhang2018subdivision}, constructions that are $C^k$ almost everywhere~\cite{sangalli2016unstructured,buchegger2016adaptively}, $G^k$ constructions~\cite{scott2013isogeometric,scott2014isogeometric, nguyen2016c,collin2016analysis,kapl2017isogeometric,kapl2017analysis,kapl2018construction} and~$C^k$ constructions with singular parameterisations~\cite{toshniwal2017smooth, toshniwal2017multi}.  Since most conventional CAD models are based on B-spline or NURBS surfaces, any isogeometric analysis suitable construction should be able to reproduce tensor-product B-splines and NURBS. For this reason we especially focus on manifold-based surfaces that can reproduce tensor-product B-splines in regular portions of the mesh. 

%

In this paper, we exploit the flexibility of the partition of unity method to devise manifold-based basis functions that can reproduce or are identical to B-splines.  The proposed techniques are introduced for the sake of clarity with the help of univariate B-splines.  Evidently, manifold-based basis functions can reproduce B-splines only on structured regions of an unstructured mesh with extraordinary vertices, i.e. non-boundary vertices with different than four attached elements. In the vicinity of extraordinary vertices the basis functions consist of a local polynomial approximant that can smoothly blend with the surrounding B-spline reproducing basis functions. The extent of the transition region depends on the size of the chosen chart domain, which consists of an $n_v$-ring of elements around each vertex. We consider several different choices for the weight functions that lead to B-spline reproducing basis functions. Especially promising are weight functions which are defined as a linear combination of B-splines defined on a grid obtained by subdividing the elements multiple times. They provide a partition of unity without normalisation and, hence, lead to polynomial manifold-based basis functions. To obtain manifold-basis functions that are identical to B-splines the local polynomial approximants have to be altered. Whereas in the original manifold constructions the local polynomial approximants are~$C^\infty$, they have to be chosen to have the same smoothness as the considered B-splines. 

On structured meshes the approximation properties of manifold-based basis functions can be inferred from the theory presented in Melenk and Babuska~\cite{Melenk1996a}. The summation of the local errors on the charts gives a global error estimate under some smoothness assumptions on the weight functions. The local error, for instance in~$L^\infty$ or~$L^2$ norms, is bounded by $h^{p+1}$, where~$h$ is the diameter of the chart domain and~$p$ is the degree of the polynomials contained in the local approximant. Each chart domain consists of $n_v$-rings of elements around a vertex so that there are on unstructured meshes multiple types of chart domains depending on the local connectivity of the control mesh. The local connectivity of the control mesh determines the type of transition function used in the manifold construction. When a control mesh is refined by quadrisecting its elements, all the newly introduced vertices are ordinary. That is, for points close to the extraordinary vertices the type of the transition function used depends on the refinement level of the mesh. Therefore, the theory presented in Melenk and Babuska~\cite{Melenk1996a} has to be extended to cover the extraordinary vertices, which we do not attempt in this paper.  

The outline of this paper is as follows. In Section~\ref{sec:review} we review the manifold basis functions as introduced in~\cite{majeedCirak:2016}. Although only univariate basis functions on polygonal control meshes are considered, no specific choices for the weight functions, local approximants and the transition functions are given so that the presented theory is applicable to the multivariate case as well. Subsequently, in Section~\ref{sec:reproduction} several specific choices first for weight functions and then for local approximants are introduced. More specifically, in Section~\ref{sec:reprodW} five choices for the weight functions are proposed, two of which yield manifold basis functions that can reproduce B-splines. In Section~\ref{sec:reprodP} it is illustrated how to choose the local approximants so that manifold-basis functions are identical to B-splines. Finally, in Section~\ref{sec:conclusions} we provide a summary and comparison of the different proposed constructions. 

%
\section{Review of manifold-based basis functions}~\label{sec:review}
%
We provide a brief informal review of univariate manifold basis functions for curves with the aim to fix ideas and notation. 
%
\subsection{Basic approach}~\label{sec:reviewBasic}
%
Given is a control polygon with the vertex coordinates~$\vec x_i \in \mathbb R^3$  which describes the manifold curve~$\Omega$. To begin with, the control polygon and the curve are assumed to be closed to sidestep the discussion of boundaries. The curve is composed of a set of~$n_c$ overlapping subdomains
\begin{equation}
	\Omega = \bigcup_{i =1}^{n_c} \Omega_i \, .
\end{equation}
Each subdomain~$\Omega_i$ is associated with a vertex~$\vec x_i$ of the control polygon in a manner yet to be described.  The subdomains~$\Omega_i$ are obtained from corresponding planar domains~\mbox{$\hat \Omega_i \in \mathbb R$} with a mapping 
%
\begin{equation}\label{eq:varphi}
    \begin{array}{rrl}
	\vec \varphi_i  \colon  &\hat \Omega_i &\rightarrow \Omega_i \\
	&\xi_i & \mapsto \vec x \, .
	\end{array}
\end{equation}
The pair consisting of $( \hat \Omega_i \, , \vec \varphi_i)$ is called a chart.  In the following we refer to~$\hat \Omega_i$ as the chart domain or simply as the chart. If it is clear from context, we neglect the mapping~$\vec \varphi_i$. Each chart domain~$\hat \Omega_i$ has its own coordinate system with the coordinates~$\xi_i \in \mathbb R$. The coordinates of points on  the intersection between two subdomains~$\Omega_i$ and~$\Omega_j$ can be transformed with  transition functions, that is, 
%
\begin{equation} \label{eq:transitionF}
    \begin{array}{rrl}
	t_{ji}  \colon  &\hat \Omega_i^j \subset \hat \Omega_i &\rightarrow \hat \Omega_j \\
	&\xi_i & \mapsto \xi_j \, 
	\end{array}
\end{equation}
defined as
\begin{equation}  \label{eq:transitionFm}
	 t_{ji} = \vec \varphi_j^{-1} \circ \vec \varphi_i \, .
\end{equation}
Here $\hat \Omega_i^j = \vec \varphi_i^{-1}(\Omega_i \cap \Omega_j)$ is the pull-back of the intersection of the two subdomains.

For constructing a smooth approximant, on each chart domain we have given a blending (or, weight) function $w_i : \hat \Omega_i \rightarrow \mathbb{R}^+_0$ with
\begin{equation}
    \mbox{supp}(w_i)\subseteq \hat \Omega_i
\end{equation}
which have to satisfy
\begin{equation}~\label{eq:wSumTo1}
    \sum_{i=1}^{n_c} w_i\circ \vec \varphi_i^{-1} \equiv 1 \quad\mbox{ on } \Omega  
\end{equation}
and have to be~$C^k$ smooth. In addition, at the chart domain boundaries,~$w_i$ and its derivatives up to $k$-th order have to be zero. On each chart domain also a local approximant $f_i : \hat \Omega_i \rightarrow \mathbb{R}$ is defined. The approximant~$f_i$ is usually expressed in a polynomial basis, like the power, Lagrangian or the B\'ezier basis. In this paper both the Lagrangian and B\'ezier basis of a fixed degree are used. However, the basis and the degree of the approximant~$f_i$ may be different on every chart.  Hence, having the local bases $\mathcal{P}_i = \left \{ p_{i}^{(j)}(\xi_i) \right \}_{j=1}^{q_p+1}$ of degree~$q_p$ on  chart $\hat \Omega_i$ gives the local approximant 
\begin{equation}~\label{eq:fi}
	f_i(\xi_i) = \sum_{j=1}^{q_p+1}p_{i}^{(j)}(\xi_i) \alpha_{i}^{(j)} := \vec{p}^\trans_i(\xi_i) \vec{\alpha}_i \,,
\end{equation}
and, in turn, the global approximant
\begin{equation}~\label{eq:globalf}
	f(\xi_i) = \sum_{l \colon \vec \varphi_i (\xi_i) \in \vec \varphi_l(\xi_l)} w_l(\xi_l)f_l(\xi_l) = \sum_{l \colon \vec \varphi_i (\xi_i) \in \vec \varphi_l(\xi_l)}  w_l(\xi_l) \vec{p}^\trans_l(\xi_l) \vec{\alpha}_l  \,,
\end{equation}	
as well as the global basis
\begin{equation}\label{eq:globalBasis}
    \mathcal{P}_{global} = \bigcup_{i=1}^{n_c} w_i \mathcal{P}_i = \{ w_i(\xi_i) p_{i}^{(j)}(\xi_i): \mbox{ with }i=1,\ldots, \,  n_c \mbox{ and }j=1,\ldots, \,  q_p+1\} \,.
\end{equation}
Note that the index~$i$ from the basis~$ p_{i}^{(j)}$ may be dropped when on each chart domain~$\hat \Omega_i$ the same basis is used, as in the present paper. 

Next, each chart domain~$\hat  \Omega_i$ and its image~$\Omega_i$ are associated with segments/elements  in the $n_v$-neighbourhood of the control polygon around the vertex~$\vec x_i$. That is, there are as many charts as vertices in the mesh. An $1$-neighbourhood of a vertex is defined as the union of elements that contain the vertex. The $n_v$-neighbourhood is defined recursively as the union of all 1-neighbourhoods of the $(n_v-1)$-neighbourhood vertices. The number of segments associated with a chart is hence~$2n_v$. In turn, each segment is present in $2n_v$ charts. See  Figure~\ref{fig:1dConstruction} for a construction with each  chart domain~$\hat  \Omega_i$  consisting of the two segments in the $1$-neighbourhood of the vertex~$\vec x_i$. 

In Section \ref{sec:reproduction} we consider the span of $\mathcal{P}_{global}$ in~\eqref{eq:globalBasis} as the \emph{analysis space} on $\Omega$. That is, the basis together with the corresponding set of coefficients~$\{ \vec \alpha_j \}$ is used for $L^2$-fitting.  However, the set of coefficients~$\{ \vec \alpha_j \}$ lack an intuitive interpretation, similar to the control vertices of splines,  so that the basis~\eqref{eq:globalBasis} is not suitable for geometric modelling. Hence, in the following we define a \emph{design space} as a suitable subspace having degrees of freedom corresponding to the vertex coordinates $\vec x_i$. 

%
\subsection{Mesh-based approach}~\label{sec:reviewMesh}
\begin{figure}
	\centering
	\includegraphics[scale=0.85]{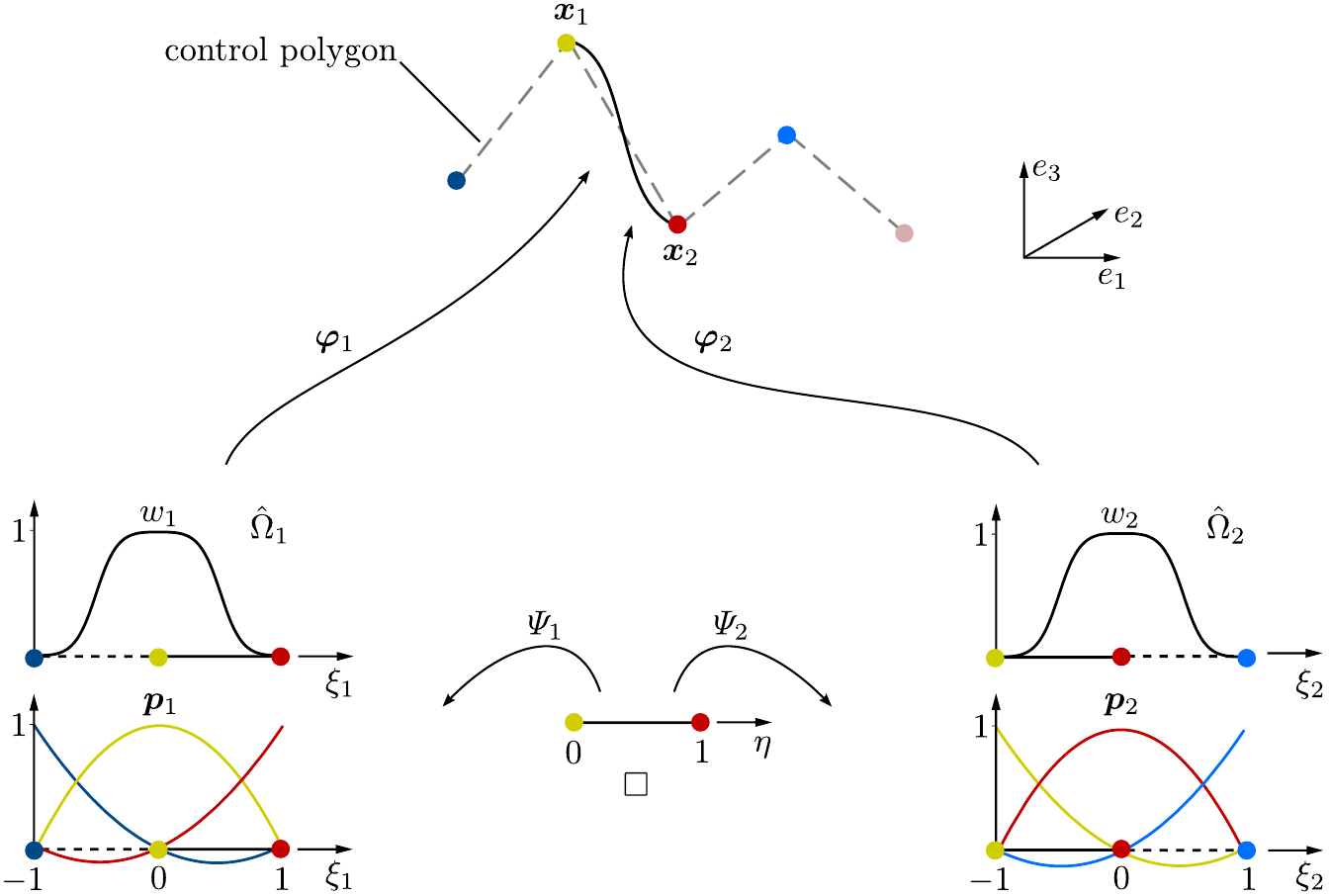}
	\caption{Construction of a univariate manifold basis. The chart domains~$\hat \Omega_1$ and~$\hat \Omega_2$ are chosen to consist out of two segments. The segment~$[\vec x_1, \, \vec x_{2}] $ on the (dashed) control polygon  is present on both chart domains.  On each chart the local basis~$\vec p_1 (\xi_1)$ and~$\vec p_2 (\xi_2)$ is a quadratic Lagrange basis and the blending functions~$w_1(\xi_1)$ and~$w_2 (\xi_2)$ are normalised B-splines. The reference element (i.e. reference finite element) is denoted with~$\Box$. The same vertex has the same colour in the four domains~$\Box$, $\hat \Omega_1$, $\hat \Omega_2$ and~$\Omega$. Note that  for the considered quadratic Lagrange basis and three vertices per chart domain the least-squares projection matrices~$\vec A_1$ and~$\vec A_2$ are identity matrices. Only choosing a constant or linear Lagrange basis on each chart domain requires a least-squares projection.}
	\label{fig:1dConstruction}
\end{figure}

%
%
On each chart domain the coefficients of the local approximant can be assigned to vertices in the $n_v$-neighbourhood, see Figure~\ref{fig:1dConstruction}. Each vertex is present on~$2 n_v+1$ charts which leads to a coupling between the coefficients of the local approximants~\eqref{eq:fi} of the involved charts. If the number of the coefficients of the local approximant is less than the number of vertices in the chart a least squares fitting has to be applied %
\begin{equation}~\label{eq:leastSquares}
	\vec{\alpha}_i  = \vec{A}_i \vec{P}_i \vec f \,,
\end{equation}
where $\vec{A}_i$ denotes the least-squares projection matrix,  \vec f is an array of the scalar vertex coefficients~$f_i \in \mathbb R$ for the entire polygon (one scalar per vertex) and $\vec{P}_i$ a gather matrix filled with ones and zeros to pick up the control vertex coefficients for a particular chart from $\vec{f}$. Hence, the global approximant~\eqref{eq:globalf} can be rewritten as 
\begin{equation} \label{eq:globalBasisLSQ}
	f(\xi_i) =  \sum_{l \colon \vec \varphi_i (\xi_i) \in \vec \varphi_l(\xi_l)}  w_l(\xi_l) \vec{p}^\trans_l (\xi_l) \vec{A}_l \vec{P}_l \vec f  \,,
\end{equation}
Finally, the manifold curve~$\Omega$ can be obtained by replacing the array of vertex scalar coefficients~$\vec f$  with the array of  given vertex coordinates~$\vec x$,  so that each map~\eqref{eq:varphi} reads 
\begin{equation}
	\vec \varphi_i   (\xi_i) =   w_i(\xi_i) \vec{p}^\trans_i (\xi_i) \vec{A}_i \vec{P}_i \vec x \, .
\end{equation}

To summarise so far, each segment on the control mesh has a unique set of corresponding segments on several planar chart domains~$\hat \Omega_i$. The introduced manifold construction ensures that the images of the set of segments from disparate planar charts are identical on the manifold~$\Omega$. To advance a more classical finite element interpretation, each segment on the manifold~$(\Omega, \,s)$ with the index~$s$ represents an element and has a corresponding reference element~\mbox{$( \Box \coloneqq [0, \, 1]  \ni \eta, \, s)$} to evaluate the element integrals\footnote{The index~$s$ for the reference element~$\Box$ is usually dropped because all of them can be assumed to have the same domain.}. The mapping of the parent element onto the manifold is composed of two maps
%
\begin{equation} \label{eq:mapBoxToChart}
    \begin{array}{rrcccc}
	\vec \varphi_i \circ \Psi_{i, s} \colon  &(\Box, \, s) &\rightarrow& (\hat \Omega_i, \, s) &\rightarrow& (\Omega, \, s) \\
	& (\eta,s) &\mapsto& \xi_i & \mapsto& \vec x 
	\end{array}
\end{equation}
with  $(\hat \Omega_i, \, s)\subset \hat \Omega_i$ and $(\Omega, \, s) \subset \Omega$ being a segment on the chart domain or manifold, respectively. 
 This implies for the field variables in a reference element~$s$, 
\begin{equation}
	f (\eta, \, s) \equiv f (\Psi_{i,s}^{-1} (  \varphi^{-1}_i  (\vec x) ) ) \, .
\end{equation}
In applications the maps~$\Psi_{i,s}$ have to be chosen carefully. Namely, the transition functions defined in~\eqref{eq:transitionFm} can be determined as, c.f.~Figure~\ref{fig:1dConstruction}, 
\begin{equation}
	t_{ji} = \Psi_{j,s} \circ \Psi_{i,s}^{-1}
\end{equation}
so that the required smoothness of~$t_{ji}$ depends on the collection of maps~$ \Psi_{i,s}$ and~$ \Psi_{j,s}$ on the respective chart domains.

The approximation of the field variables with~\eqref{eq:globalBasisLSQ} leads for the element~$s$ to the following definition of finite element basis functions~$\vec N(\eta, \, s)$: 
\begin{equation}~\label{eq:manifoldBasis}
	f(\eta , \, s) = \underbrace{   \left(\sum_{j \colon \vec \varphi_i (\xi_i) \in \vec \varphi_j(\xi_j)}  w_i(\xi_i)\vec{p}^\trans_i(\xi_i) \vec{A}_i\vec{P}_i\right)}_{\vec{N}^\trans (\eta, \, s)} \vec{f} \quad \text{with} \quad \xi_i = \Psi_{i,s}(\eta ) \,.
\end{equation}	
In the following we denote the basis functions with~$\vec{N}(\eta)$ and the mapping from the reference element to the chart with~$\Psi_s$. This notation is not precise when charts have different geometries and number of vertices. 

It is clear that the smoothness of the basis functions $\vec{N}(\eta)$ depends on the smoothness of blending functions $w_i(\xi_i)$, local basis $\vec{p}_i(\xi_i)$ as well as the mappings $\Psi_i(\eta)$. For instance, in the two-dimensional $C^k$ continuous construction introduced in~\cite{majeedCirak:2016}, the blending functions $w_i$ are chosen to be (normalised) B-splines of degree $k+1$, the local basis~$\vec p_i$ are chosen to be a polynomial basis and $\Psi_i$ are conformal maps. Furthermore, each chart domain consists of the elements in the 1-neigborhood of the corresponding vertex.   Figure~\ref{fig:1dConstruction} illustrates this construction for $C^2$ continuous basis functions in the univariate case.  
%

%
\section{Reproduction of B-splines \label{sec:reproduction}}
%
We consider again the global basis~\eqref{eq:globalBasis} for analysis, repeated here for convenience,
\begin{equation*}
    \mathcal{P}_{global} = \bigcup_{i=1}^{n_c} w_i \mathcal{P}_i = \{ w_i(\xi_i) p_{i}^{(j)}(\xi_i): \mbox{ with }i=1,\ldots, \, n_c \mbox{ and }j=1,\ldots, \, q_p+1 \} \,,
\end{equation*}
and discuss how to choose the blending functions $w_i$ and the local basis~$ p_{i}^{(j)}$ so that B-splines are a subset of the basis  $\mathcal{P}_{global} $. To reproduce B-splines it is sufficient to choose  either the blending functions $w_i$ or the local basis $\vec{p}_i$ suitably. In the following we introduce several choices for the blending functions and local basis and comment on their extendability to the bivariate case. Evidently, B-splines are defined on a  structured mesh  so that manifold-based basis functions will only reproduce B-splines on the parts of the mesh with no extraordinary vertices. 

In the univariate case, the parameter domain of the manifold curve~$\Omega$ can be assumed to be one single finite interval~$\hat \Omega$. Due to the choice of the single finite interval the transition functions~$t_{ij}$  are identity maps and~$\Psi_{s}$ are affine maps, both are omitted in the following.  Without loss of generality, the parameter domain is uniformly partitioned with~$n$ inner nodes with the coordinates~$\hat \xi_j =j$.  Moreover, we use the notation $\mathbb{P}^p$ for polynomials of degree $p$, $\mathcal{S}^{p,r}$ for B-splines of degree $p$ and continuity $C^r$ and define the space 
\begin{equation*}
    \mathcal{W} = \mbox{span}\{ w_i: \; i=1,\ldots, \, n_c \} \,
\end{equation*}
as a function space on~$\hat \Omega$.

%
\subsection{Blending function choices \label{sec:reprodW}}
%
We fix the local basis~$ p_{i}^{(j)}$ to be polynomials of some prescribed degree and study how to choose the  blending functions~$w_i$ to obtain B-spline reproducing manifold-based basis functions. As the B-splines form a local, non-negative partition of unity, they are used as blending functions. Specifically, the blending functions will be chosen either as
\begin{itemize}
 \item[--]standard B-splines of maximum smoothness,
 \item[--] rational B-spline functions, or
 \item[--] linear combinations of B-splines.
\end{itemize}
We compare the different approaches in terms of maximum number of overlapping charts at any point of the domain, expected approximation order as well as smoothness properties.

%
\subsubsection{Piecewise linear $C^0$ continuous blending functions \label{subsec:linear-weights}}
\begin{figure}[]
    \centering
 \includegraphics[scale=0.9]{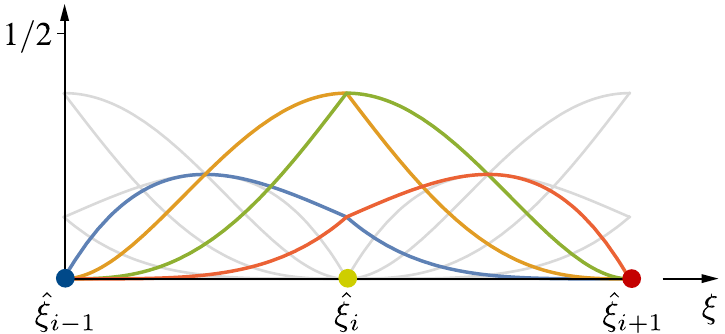}
 \caption{Global basis for piecewise linear blending function with $q_w=1$ and local cubic B\'ezier basis with $q_p=3$  restricted to the chart domain~$\hat \Omega_i = \left [ \hat \xi_{i-1}, \, \hat \xi_{i+1} \right ]$.}
 \label{fig:C0weight-basis}
\end{figure}

\begin{figure}[]
    \centering
 \includegraphics[scale=0.9]{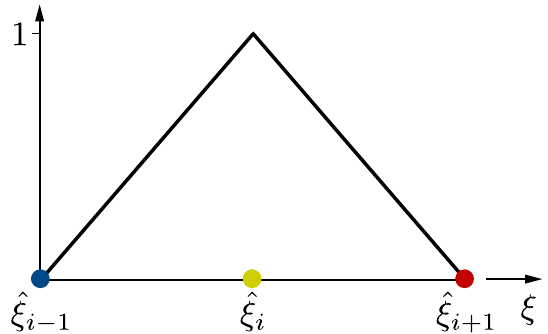} \hfill
 \includegraphics[scale=0.9]{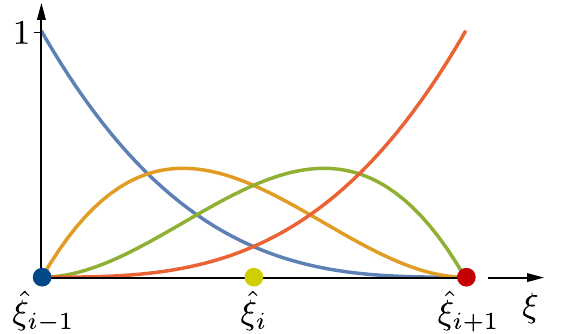}
 \caption{Piecewise linear blending function with $q_w=1$ (left) and a local cubic B\'ezier basis with $q_p=3$  (right) used in computing the basis in Figure~\ref{fig:C0weight-basis}.}
 \label{fig:hat-Bezier}
\end{figure}

%
In case of  linear B-spline blending functions of degree $q_w=1$ each chart domain contains three knots. This means that any point on the parameter domain is present on two different chart domains. Having $C^0$ hat functions as blending functions and polynomials of degree $q_p$ as local functions on every chart, we reproduce continuous, piecewise polynomials of degree $q_p+1$, i.e., $C^0$ B-splines. Hence, $\mathcal{W} = \mathcal{S}^{1,0}$ and $\mbox{span}(\mathcal{P}_{global}) = \mathcal{S}^{q_p+1,0}$.  In Figure~\ref{fig:C0weight-basis} the manifold-based basis obtained with linear B-spline blending functions and a cubic B\'ezier local basis is shown. The corresponding hat function and the cubic B\'ezier basis are depicted in Figure~\ref{fig:hat-Bezier}.

Extending the construction to surfaces, we obtain functions that are tensor-product polynomials within every quadrilateral element and $C^0$ continuous across every edge. Hence, we can only reproduce $C^0$ continuous basis functions on unstructured quadrilateral B\'ezier meshes. 

%
\subsubsection{Higher order $C^{p-1}$  continuous B-spline blending functions \label{subsec:smooth-spline-weights}}
\begin{figure}[]
    \centering
 \includegraphics[scale=0.875]{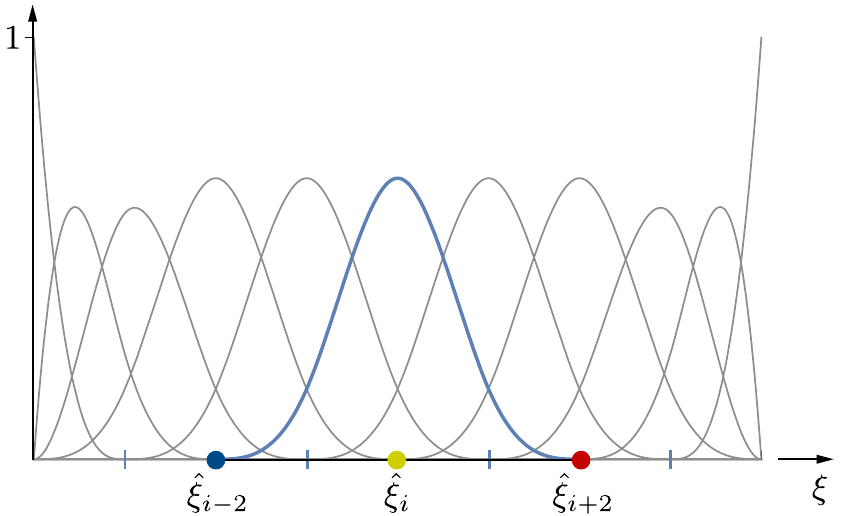}\hspace{.1\textwidth}
 \includegraphics[scale=0.875]{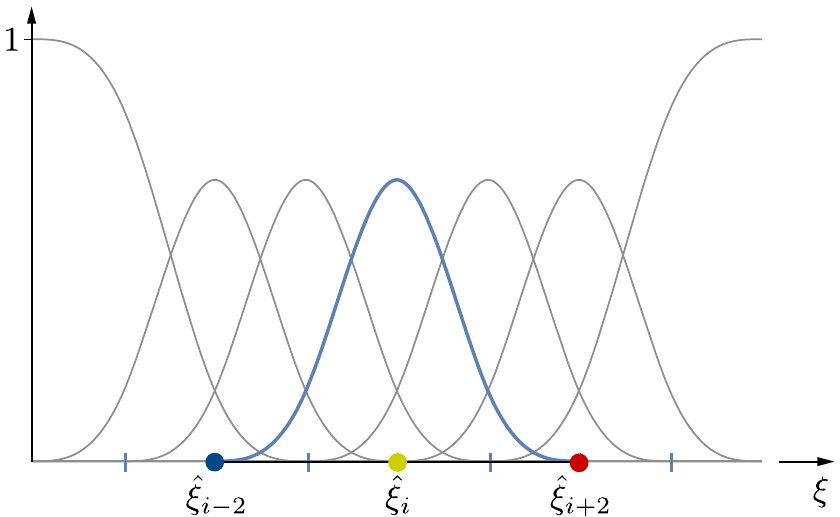}
 \caption{Spline blending functions with~$q_w=3$ without boundary correction (top), with boundary correction (bottom).}
 \label{fig:Bsplines-boundary}
\end{figure}

%
The generalisation of the  linear B-spline blending functions to the higher order B-splines with $q_w>1$ is straightforward. As shown in Figure~\ref{fig:Bsplines-boundary} the support of each basis function on the parameter domain is defined as a chart. Therefore, for  B-splines of degree $q_w$ and smoothness $q_w-1$ in any point of the domain $q_w+1$ charts overlap. 

Taking splines of degree $q_w >1$ as blending functions, together with polynomials of degree $q_p$ as local basis, we reproduce B-splines of degree $q_w+q_p$ and smoothness $q_w-1$. Hence, we have $\mathcal{W} = \mathcal{S}^{q_w, \, q_w-1}$. 
Let $\mathcal{W}_0$ be the subspace of $\mathcal{W}$ without global polynomials, i.e.,
\begin{equation*}
\mathcal{W}_0 = \mathcal{W} / \mathbb{P}^{q_w},
\end{equation*}
then we have 
\begin{equation*}
\mbox{span}(\mathcal{P}_{global}) = \mathbb{P}^{q_w+q_p} \oplus \xi^{q_p} \mathcal{W}_0 \oplus \ldots \oplus \xi \mathcal{W}_0 \oplus \mathcal{W}_0 = \mathcal{S}^{q_w+q_p, q_w-1}.
\end{equation*}
The dimension of $\xi^k \mathcal{W}_0$ is independent of $k$ and is equal to the number $n$ of inner knots of the spline space $\mathcal{W} = \mathcal{S}^{q_w,q_w-1}$, i.e., $\mbox{dim}(\mathcal{W}) = n+q_w+1$, $\mbox{dim}(\mathcal{W}_0) = n$ and 
\begin{equation*}
\mbox{dim}(\mbox{span}(\mathcal{P}_{global})) = (q_w+q_p+1) + n\cdot (q_p+1).
\end{equation*} 
As the span of $\mathcal{P}_{global}$ contains splines of degree $q_w+q_p$, we can expect an approximation order of $O(h^{q_w+q_p+1})$ in $L^2$. 
Note that since the space of blending functions contains polynomials, 
the functions in $\mathcal{P}_{global}$ are linearly dependent for $q_w >0$, as we then have $|\mathcal{P}_{global}| = (q_p+1)\cdot (n+q_w+1) > (q_w+q_p+1) + n\cdot (q_p+1)$.

The treatment of the boundary is not straightforward. In Figure \ref{fig:Bsplines-boundary} two different options for choosing cubic blending functions with $q_w=3$ is presented. The treatment of boundaries becomes relevant when extending the construction to surfaces. This leads to surfaces that are $C^{q_w-1}$ if the mesh is regular. The construction on the top contains all splines, but leads to $C^0$ smooth surfaces at the extraordinary vertices. The construction on the bottom reproduces only a subspace of all B-splines, but generates surfaces that are $C^2$ smooth everywhere. However, as the overlap between the charts is large, the construction becomes cumbersome. Especially, in the bivariate case there can be several extraordinary vertices within one chart domain which can render their smooth parametrisation challenging.

For this reason, in the following, we consider blending functions that have a small support, but generate a smooth basis.
%
\subsubsection{Rational B-spline blending functions}\label{subsec:rational-weights}
\begin{figure}[]
    \centering
 \includegraphics[scale=0.9]{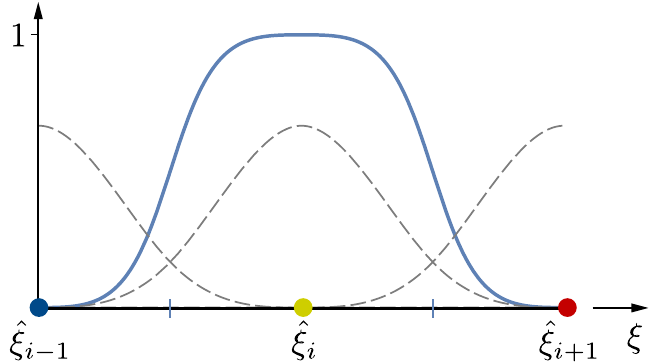}
 \caption{Blending function~$w_i$  (solid line) as the sum  of a selection of cubic B-splines (dashed lines).}
 \label{fig:selected-Bsplines}
\end{figure}

%
To circumvent the difficulties that arise from using chart domains with a large number of overlaps, we construct blending functions that lead to only two overlapping chart domains in any point of the domain. That is, we consider blending functions with a overlaps similar to the hat functions in Figure \ref{fig:hat-Bezier} (left), but possess the smoothness of higher order B-spline blending functions as in Figure \ref{fig:Bsplines-boundary} (bottom). To reproduce this behaviour, we first select a linear combination of the B-splines in Figure \ref{fig:Bsplines-boundary} as blending functions, so that that the supports of no more than two blending functions overlap at the same time. In addition the B-splines to be used as blending functions are chosen from a suitably scaled coordinate system. See Figure \ref{fig:selected-Bsplines} for a construction with cubic B-splines with $q_w=3$ defined on a coordinate axis scaled by a factor~$2$. Here, the blending functions are defined as 
\begin{equation}\label{eq:rationalweight}
 w_i (\xi) = \frac{B_{2i}^3 (2\xi)}{ B_{2(i-1)}^3  (2 \xi) + B_{2i}^3  (2 \xi) + B_{2(i+1)}^3  (2 \xi)}, 
\end{equation}
where $B_k^3$ are the B-spline basis functions. Note that every other B-spline is chosen as a blending function and their sum does not add up to one. To obtain a partition of unity, all functions are divided by their sum. The resulting blending functions are then piecewise rational, as depicted in Figure \ref{fig:selected-Bsplines}. The resulting manifold-based basis functions are also rational and their numerical integration may need more quadrature points than polynomial basis functions of similar order and smoothness. 

%
\subsubsection{Linear combinations of B-splines as blending functions}\label{subsec:lc-weights}
\begin{figure}[]
    \centering
 \includegraphics[scale=0.9]{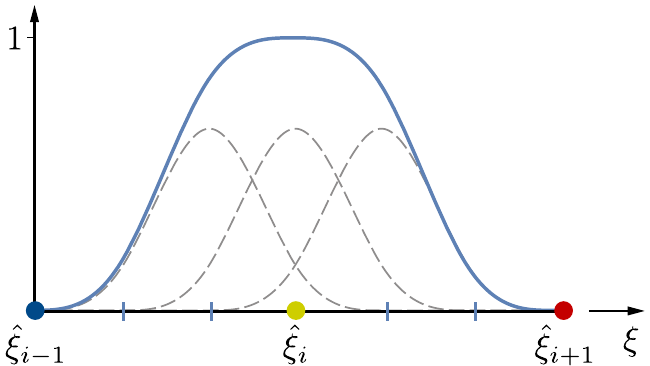}
 \caption{Blending function~$w_i$  (solid line) composed out of three consecutive cubic B-splines (dashed lines) and the mask $(1, \, 1, \, 1)$.}
 \label{fig:sum-Bsplines}
\end{figure}

\begin{figure}[]
    \centering
 \includegraphics[scale=0.9]{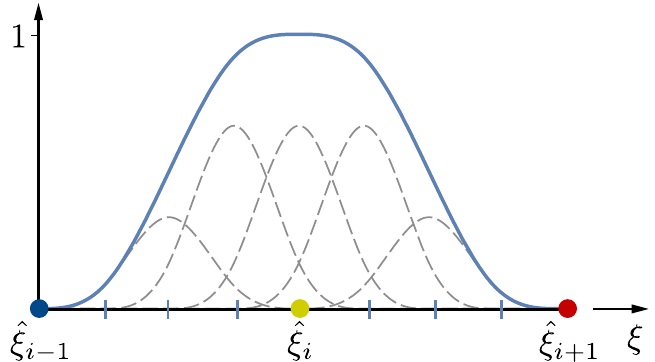}
 \caption{Blending function~$w_i^*$  (solid line) composed out of five consecutive cubic B-splines (dashed lines) and the mask $(\frac{1}{2}, \, 1, \, 1, \, 1,\frac{1}{2})$.}
 \label{fig:lc-Bsplines}
\end{figure}

%
One can take linear combinations of consecutive B-splines as blending functions to obtain polynomial blending functions. We consider only cubic B-splines with~$q_w=3$ and express the coefficients in the linear combination as masks $(m_0, \, m_1, \,\dotsc, \, m_k)$.   To begin with, the B-splines to be used as blending functions are defined on a coordinate axis scaled by a factor~$3$ and the blending function mask is  $(1, \, 1, \, 1)$, see Figure \ref{fig:sum-Bsplines}. That is, the blending functions are obtained as
\begin{equation}\label{eq:weight111}
 	w_i (\xi)= B_{3i}^3 (3\xi)+ B_{3i+1}^3 (3\xi) + B_{3i+2}^3 (3\xi)  \, .
\end{equation}
Alternatively, it is possible to use B-splines defined on a coordinate axis scaled by a factor~$4$ with a mask $(\frac{1}{2}, \, 1, \, 1, \, 1, \, \frac{1}{2})$ so that
\begin{equation}\label{eq:weight1112}
	 w_i^* (\xi)= \frac{1}{2}B_{4i}^3 (4\xi)+ B_{4i+1}^3 (4\xi) + B_{4i+2}^3 (4\xi)+ B_{4i+3}^3 (4\xi)+ \frac{1}{2}B_{4i+4}^3 (4\xi),
\end{equation}
This choice is illustrated in Figure~\ref{fig:lc-Bsplines}.  As is indicated in the Figures~\ref{fig:sum-Bsplines} and~\ref{fig:lc-Bsplines} within a chart the blending function~$w(\xi)$ and~$w^*(\xi)$ have five  and seven breaking points, respectively.  


\begin{figure}[]
    \centering
 \includegraphics[scale=0.9]{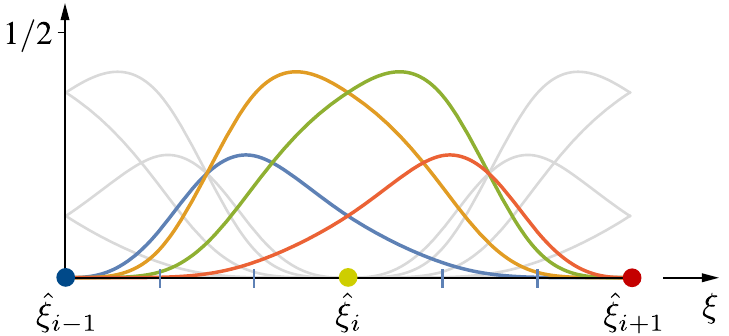}
 \caption{Global basis for a B-spline blending function composed according to~\eqref{eq:weight111} and a local cubic B\'ezier basis with~$q_p = 3$ restricted to one chart.}
 \label{fig:basis-sum-Bsplines}
\end{figure}

We can now compute the local contribution of the global basis on one chart. This is depicted in Figure~\ref{fig:basis-sum-Bsplines}. The basis functions using weights as in~\eqref{eq:rationalweight} or~\eqref{eq:weight1112} are visually indistinguishable and are omitted here.

%
\subsubsection{Comparison of different blending function choices}
%
To summarise, the constructions using rational functions or sums of B-splines generate charts that have only a one ring overlap. The same is true for piecewise linear hat functions as blending functions. Moreover, the dimension of the global function space is the same in these three cases. However, for piecewise linears, the functions are only $C^0$, whereas for the other three approaches the smoothness is $C^{q_p-1}$. In the following we compare the approximation power of the respective approaches. 

\begin{figure}[]
    \centering
 \includegraphics[width=.925\textwidth]{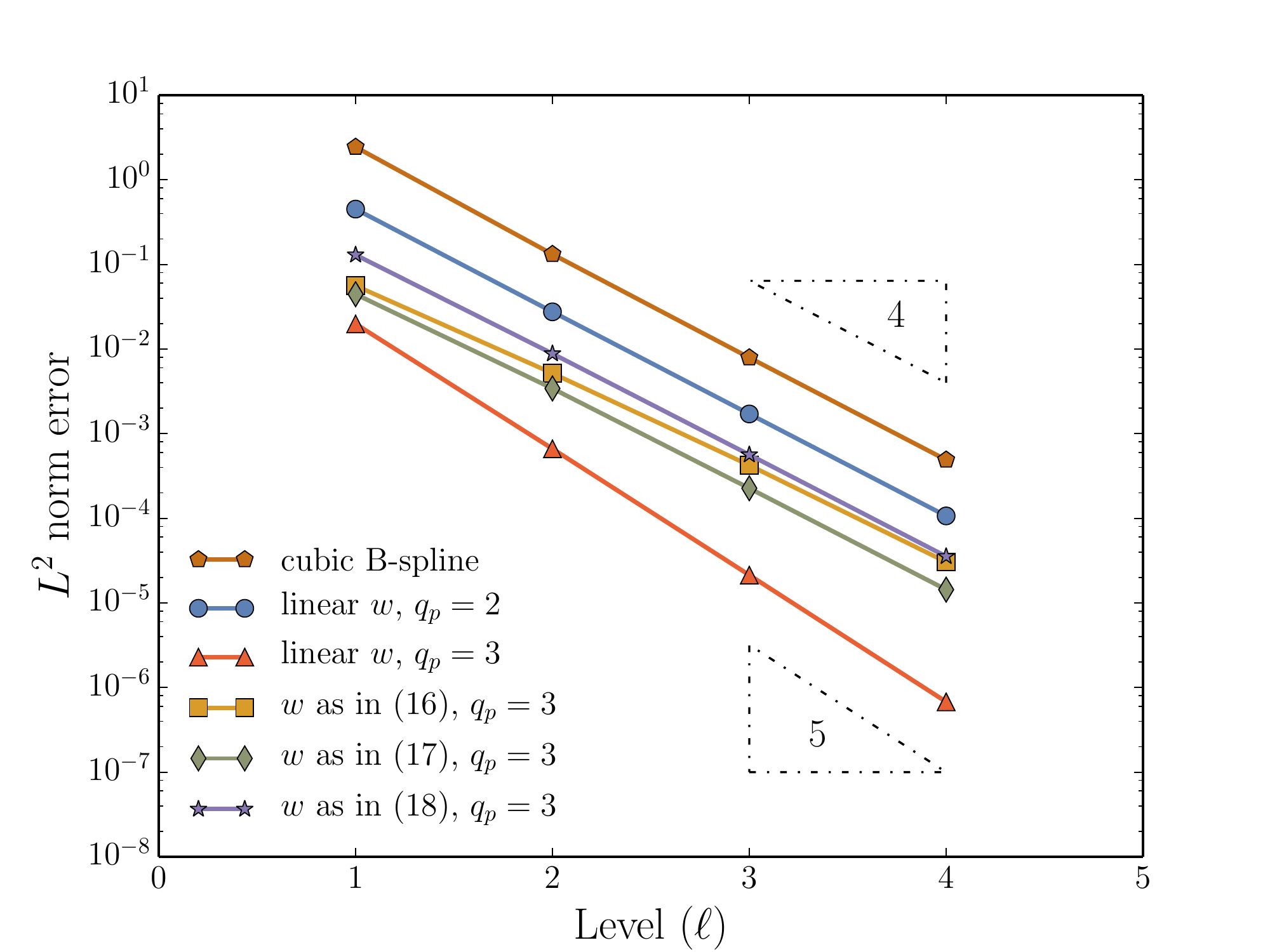}
 \caption{$L^2$ approximation error plots for various choices of blending functions. Rates for $h^4$ and $h^5$ are included for comparison.}
 \label{fig:approx-error}
\end{figure}

In Figure~\ref{fig:approx-error} we show log-error plots when performing $L^2$-fitting onto a given function. In this example we considered a sine function over the unit interval. Here, the mesh size satisfies $h=1/2^{\ell+2}$, where we used levels $\ell=1,\ldots,4$. We compare linear blending functions with local polynomials of degree $q_p=2$ (blue line) and $q_p=3$ (red line), as in Subsection~\ref{subsec:linear-weights}. The former has a theoretical convergence rate of $O(h^4)$ in $L^2$, while the latter has a theoretical rate of $O(h^5)$. Both discretizations are $C^0$ only. We moreover compare rational blending functions as in~\eqref{eq:rationalweight} (yellow line), linear combinations of splines as blending functions as in~\eqref{eq:weight111} (green line) or in~\eqref{eq:weight1112} (purple line) of degree $q_w=3$ and local polynomials of degree $q_p=3$. In all three cases, the expected convergence rate is $O(h^4)$.

We compare all constructions with uniform cubic B-splines of mesh size $h$ (orange line). They can be interpreted as a manifold construction as in Subsection~\ref{subsec:smooth-spline-weights} with $q_w=3$ and $q_p=0$. This construction yields the highest error. Similarly high errors are observed for the basis from Subsection~\ref{subsec:linear-weights} with $q_p=2$ (resulting in piecewise cubics) and the lowest error for the same construction with $q_p=3$ (resulting in piecewise quartics). The bases constructed in Subsection~\ref{subsec:smooth-spline-weights} produce rates depending on the polynomial degree of the resulting splines. As all splines of a given degree and varying smoothness converge similarly (here e.g. $q_w=2$, $q_p=1$), we have omitted this case in Figure~\ref{fig:approx-error}. Note that the weight functions with boundary correction in Figure~\ref{fig:Bsplines-boundary} (bottom) will not converge optimally without increasing the degree of local functions close to the boundary.

When comparing the constructions from Subsections~\ref{subsec:rational-weights} and~\ref{subsec:lc-weights}, it turns out that all three converge with optimal rates of order $h^4$ and with significantly smaller constant when compared to piecewise polynomials of degree $3$. This means that, when fitting onto a smooth, univariate function, manifold constructions yield a better approximation than standard B-splines, even though the manifold constructions do not reproduce B-splines. Among the manifold constructions, the non-rational variant~\eqref{eq:weight111} from Subsection~\ref{subsec:lc-weights} seems to be the faster. 

It is reasonable to assume that a B-spline compatible manifold construction further improves the approximation properties. Therefore we introduce constructions with modified local functions that reproduce B-splines.

%
\subsection{Local approximants  \label{sec:reprodP}}
%
We discuss next how to choose the local basis~$p_i^{(j)}$ so that the manifold-based construction reproduces B-splines of maximum smoothness. Here, the blending functions~$w_i$ have only to satisfy the partition of unity property. Hence, any one of the blending functions introduced in Section~\ref{sec:reprodW} can be used.  To avoid the complications arising from chart domains~$\hat \Omega_i$  with large number of overlaps, we consider only blending functions which lead to two overlapping charts. In addition, for the sake of concreteness we focus in the following on cubic B-splines and note that the proposed construction carries  over to arbitrary degree. 

\begin{figure}[]
    \centering
	 \includegraphics[scale=0.9]{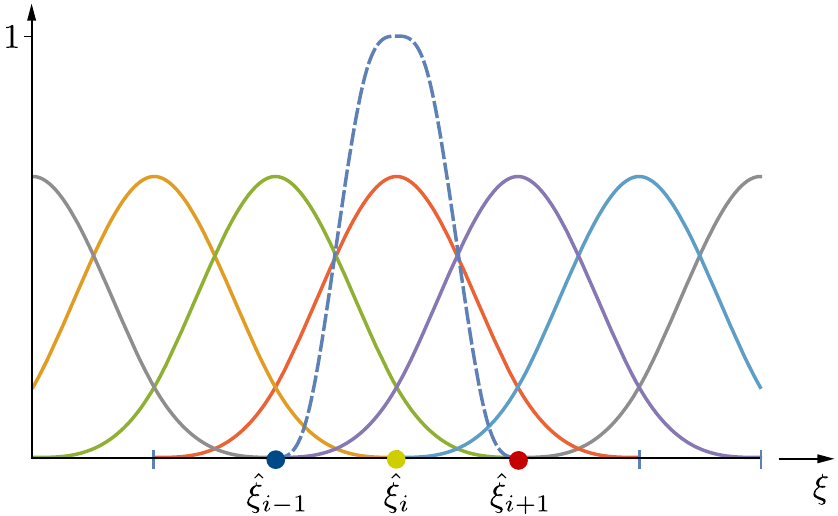}
	 \includegraphics[scale=0.9]{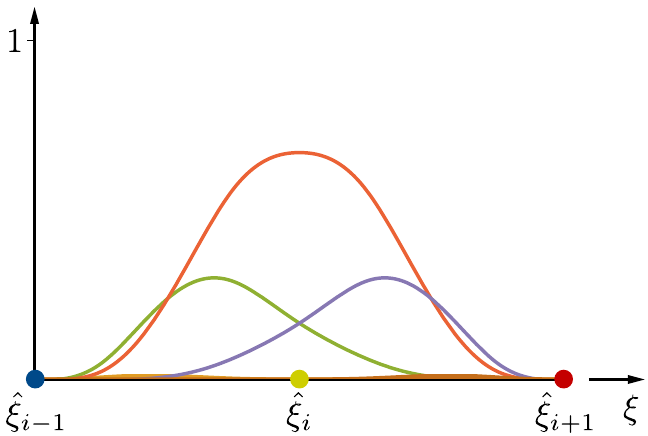}
 	\caption{Cubic spline basis functions and blending function (top). The weighted cubic spline basis functions on one chart to be reproduced with manifold-based basis functions (bottom). \label{fig:bSplineAndWeight}}
\end{figure}

The global manifold-based approximant~\eqref{eq:globalf} on a parameter domain consisting of a single finite interval is given by
\begin{equation} \label{eq:approxPc}
	f(\xi) =  \sum_{i}  w_i(\xi)  f_i (\xi)  = {  \sum_{i} } w_i(\xi)  \left (  \sum_j p_i^{(j)}(\xi) \alpha_i^{(j)} \right ) \, .
\end{equation}
It is required that this approximant is equal to a B-spline over all or some of the chart domains~$\hat \Omega_i \equiv \supp \omega_i$. The cubic B-spline approximant is defined as 
\begin{equation} \label{eq:bsplinePc}
	f^{B} (\xi) = \sum_k B_k^3 (\xi) \beta_k \, ,
\end{equation}
where~$\beta_i$ are the control point coefficients. This approximant can be expressed as a weighted sum of chart domain contributions by multiplication with the partition of unity function, that is, 
\begin{equation}  \label{eq:bsplineSplitPc}
	f^{B} (\xi) =  \underbrace{ \sum_i w_i (\xi)}_{\equiv 1}  \sum_k B_k^3(\xi) \beta_k = \sum_i w_i(\xi) f_i^B(\xi)
\end{equation}
The support of each of the terms~$w_i(\xi) f_i^B(\xi)$ is strictly restricted to one chart domain, see Figure~\ref{fig:bSplineAndWeight}. Note that the local basis in Figure~\ref{fig:bSplineAndWeight} (bottom) consists of five functions, which are scaled differently due to the multiplication with the blending function. 
Term by term matching of the manifold-based~\eqref{eq:approxPc} and the weighted B-spline approximants~\eqref{eq:bsplineSplitPc} requires on every chart domain 
\begin{equation} \label{eq:bsplineEquaPc}
	 \sum_j p_i^{(j)} \alpha_i^{(j)} =  f_i^B(\xi) \quad \text {for } \xi \in \hat \Omega_i = [ \hat \xi_{i-1}, \, \hat \xi_{i+1} ]  \, .
\end{equation} 
This equation yields a set of equations for determining the coefficients~$ \alpha_i^{(j)}$ in dependence of the coefficients~$\beta_k$. The cubic B-spline~$ B_k^3(\xi)$ has one knot with~$C^{k}$, $k\le2$, at the centre of the chart domain~$\xi=\hat \xi_i$ so that~$p_i^{(j)}$ has to consist out of two pieces. Indeed, it is sufficient to consider in each half~$[\hat \xi_{i-1}, \, 0]$  and~$[0, \, \xi_{i+1}]$ of the chart domain~$\hat \Omega_i$ a separate polynomial approximant. Choosing in each half a B\'ezier basis the coefficients~$\alpha_i$ can simply be obtained by B\'ezier extraction as a linear combination of the B-spline coefficients~$\beta_k$. Or more generally, the coefficients~$\alpha_i$  are obtained by solving a small linear system of equations obtained by collocating~\eqref{eq:bsplineEquaPc} at four distinct points (for a cubic B-spline) within the segment. Since~$ f_i^B(\xi) $ is not known~\eqref{eq:bsplineEquaPc} has to be considered for each of the four non-zero B-spline basis functions individually. Similar to~\eqref{eq:leastSquares}, this gives a relation between the two sets of coefficients expressed as
 \begin{equation}
 	\vec \alpha_i = \vec A \vec P_i \vec \beta  \, , 
 \end{equation}
Note that the projection matrix~$\vec A$ depends on the specific local basis chosen and is here the same on all the chart domains. Introducing the obtained coefficients into~\eqref{eq:approxPc} yields the manifold-based basis functions, which are by design the same as the B-spline basis functions.

In the bivariate case the B-spline approximant~\eqref{eq:bsplinePc} is only available on parts of the mesh with a tensor-product structure. In the vicinity of extraordinary vertices there is no representation as in~\eqref{eq:bsplineEquaPc}. In such regions, as in the original manifold construction introduced in Section~\ref{sec:reviewMesh} a continuous polynomial approximant has to be fitted  to the control mesh coefficients. The manifold construction ensures that the global approximant has the desired smoothness properties.

%
\section{Conclusions \label{sec:conclusions}}
%
We developed new manifold-based b-spline basis functions by using the manifold-based surface construction techniques from geometric modelling. As illustrated the manifold-based surface construction techniques can be understood as the extension of the partition of unity method to manifolds. Specific choices for the blending functions and local approximants yield B-splines on structured control meshes. Due to the flexibility of the partition of unity method several such choices are possible. We introduced in total five different choices for the blending functions two of which reproduce B-splines. In addition, we introduced one choice for the local approximant that leads to B-splines. 

In Table~\ref{tab:concl} the properties of the manifold-based basis functions obtained from each of the  six different choices are listed. For finite elements polynomial basis functions are to be preferred because they usually require less quadrature points to integrate. The number of breaking points within a finite element gives out of how many smoothly attached pieces a basis function consists.  For efficient numerical integration the breaking points of the basis function have to be considered so that constructions with fewer breaking points are to be preferred. In the multivariate case, on unstructured meshes only constructions which require only one-ring of elements around each vertex as a chart domain are viable.  If the chart domain consists out of more than one ring of elements, there can be several extraordinary vertices in a chart which makes their parametrisation challenging. In the regular setting the approximation order of the introduced constructions can be inferred from Melenk and Babuska~\cite{Melenk1996a}. In Table~\ref{tab:concl} the higher order convergence of the first two constructions is remarkable. The first construction yields however only~$C^0$ basis functions and the second construction requires charts with  several rings of elements. Overall, the most promising constructions for finite elements appear to be the blending functions assembled from B-splines introduced in Section~\ref{subsec:lc-weights} and the local approximant introduced in Section~\ref{sec:reprodP}. We note that the smoothness of the two resulting basis functions is~$C^k$ with~\mbox{$k=\min(q_p, \, q_w)$.} In closing, we note that the mathematical and numerical study of the introduced constructions on unstructured meshes provides a promising area for future research. 

\setlength\tabcolsep{0.32em}
  \begin{table}[tb]
\centering
\begin{tabular}{l |c | c |c | c| c}
 & rational / & B-spline  &  chart size & breaking  & approx.  \\
  &polynomial & reproducing  &     &  points & order  \\ 
  \hline 
$w$  in Sect. \ref{subsec:linear-weights} & p &  \checkmark  &  one-ring  &  0  &  $h^{q_p+2}$ \\ [0.1em] 
$w$  in Sect. \ref{subsec:smooth-spline-weights} & p & \checkmark &   $(q_w+1)/2$-ring & 0  &  $h^{q_w+q_p+1}$ \\  [0.1em] 
$w$  in Sect. \ref{subsec:rational-weights} & r & &  one-ring & 1  &   $h^{q_p+1}$ \\ [0.1em]  
$w$  in Sect. \ref{subsec:lc-weights} & p  &  &     one-ring & 2  & $h^{q_p+1}$  \\  [0.1em] 
$w^*$  in Sect. \ref{subsec:lc-weights} & p &  &   one-ring &  3 & $h^{q_p+1}$  \\  [0.1em] 
$p$ as in Sect. \ref{sec:reprodP} & p &  \checkmark &   one-ring &$ \ge 1$  &  $h^{q_p+1}$ \\  [0.1em] 
\end{tabular}
\caption{Comparison of the properties of the introduced manifold-based basis functions on structured control meshes. \label{tab:concl}}
\end{table}


\section*{Acknowledgement}
The authors would like to acknowledge the kind hospitality of the Erwin Schr\"odinger International Institute for Mathematics and Physics (ESI), where some of this research was carried out as part of the thematic programme {\it Numerical Analysis of Complex PDE Models in the Sciences}.

\bibliography{manifoldS.bib}
\bibliographystyle{styles/spmpsci}

\end{document}